\documentclass[a4paper]{article}
\usepackage[12pt]{extsizes}

\usepackage[T2A]{fontenc}
\usepackage[utf8]{inputenc}
\usepackage[russian]{babel}

\usepackage{hyphenat}
\usepackage{amssymb,amsmath,amsthm} 
\usepackage{libertinus} 
\usepackage{mathrsfs} 
\usepackage{xcolor} 
\usepackage{todonotes} 

\usepackage{geometry} 
\usepackage{hyperref} 
\hypersetup{
    bookmarks=true,
    unicode=true,
    colorlinks=true, 
    linkcolor=blue, 
    citecolor=red, 
    urlcolor =purple, 
    linkbordercolor=blue,
    citebordercolor=red,
    urlbordercolor =purple
 }

\usepackage{tikz} 
\usepackage{graphicx}
\newtheoremstyle{break}
  {}
  {}
  {\itshape}
  {}
  {\bfseries \vspace{.3em}}
  {}
  {\newline}
  {}

\theoremstyle{break}
\newtheorem{prop}{Утверждение}[section]
\newtheorem{lemm}[prop]{Лемма}
\newtheorem{theo}[prop]{Теорема}
\newtheorem{coro}[prop]{Следствие}

\theoremstyle{remark}
\newtheorem{rem}[prop]{Замечание}

\theoremstyle{definition}
\newtheorem{defi}[prop]{Определение}

\renewcommand{\_}{\kern-1.5pt\textunderscore\kern-1.5pt}

\newenvironment{myproof}
{\vspace{-1.2em} \begin{proof}}
{\end{proof}}

\renewcommand{\le}{\leqslant}

\allowdisplaybreaks 
\title{\sffamily Наследственная неразрешимость фрагментов некоторых элементарных теорий}

\author{\sffamily \bfseries Владимир Евгеньевич Карпов}

\date{\normalsize
{\sffamily Московский физико-технический институт,\\
Московская облаcть, г.\ Долгопрудный}\\[1em]
{\itshape Научный руководитель: Станислав Олегович Сперанский}
}

\begin{document}

\setcounter{page}{2}
\maketitle

\begin{abstract}
Как известно, при наличии интерпретации класса структур $\mathcal{K}_1$ в классе структур $\mathcal{K}_2$ наследственная неразрешимость (фрагмента) теории $\mathcal{K}_1$ влечёт наследственную неразрешимость (подходящего фрагмента) теории $\mathcal{K}_2$. В настоящей работе строится $\Sigma_1$-ин\-тер\-пре\-та\-ция класса всех конечных двудольных графов в классе всех пар эквивалентностей на общем конечном носителе; отсюда получается наследственная неразрешимость $\Sigma_2$-те\-о\-рии второго класса. Далее, строится $\Sigma_1$-ин\-тер\-пре\-та\-ция класса всех пар эквивалентностей на общем конечном носителе
в классе всех пар, состоящих из линейного порядка и эквивалентности на общем конечном носителе;
это даёт наследственную неразрешимость $\Sigma_2$-те\-о\-рии второго класса. Соответствующие результаты являются в известном смысле оптимальными, поскольку $\Pi_2$-те\-о\-рии рассматриваемых классов разрешимы.

\medskip
{\footnotesize \emph{Ключевые слова:}
неразрешимость,
элементарные теории,
префиксные фрагменты}
\end{abstract}

\newpage
\tableofcontents


\newpage
\section{Введение}

После формализации понятия алгоритма (см.\ \cite{Turing}) у учёных появилась возможность доказывать не~только существования алгоритмов, решающих те или иные задачи, но и отсутствие алгоритмов решения конкретных задач. В частности, известным примером неразршешимой задачи является проблема остановки для машин Тьюринга.

Вскоре после этого возникли вопросы о разрешимости элементарных теорий различных естественных классов структур (графов, решёток, групп, колец и так далее). Здесь можно считать, что задача ставится следующим образом: имеется класс структур $\mathcal{K}$, чьей \emph{теорией} называется множество всех первопорядковых предложений, истинных во всех структурах из $\mathcal{K}$; нас интересует, существует ли алгоритм, который по произвольному предложению определяет, принадлежит ли оно теории $\mathcal{K}$. Довольно быстро стало понятно, что неразрешимых теорий намного больше, чем разрешимых. Обзор известных результатов в этой области можно найти в 
\cite{ELTT}.

Для переноса результатов о разрешимости или неразрешимости с одних теорий на другие используется \emph{метод интерпретаций}, или \emph{метод относительной элементарной определимости}. Как известно, при наличии интерпретации класса структур $\mathcal{K}_1$ в классе структур $\mathcal{K}_2$ наследственная неразрешимость теории $\mathcal{K}_1$ влечёт наследственную неразрешимость теории $\mathcal{K}_2$. В частности, данный подход позволяет строить иерархии неразрешимых теорий. Исторически, в качестве исходной наследственно неразрешимой теории бралась теория конечных простых графов или теория двух эквивалентностей на общем конечном носителе.

В статье \cite{Nies} было предложено сконцентрироваться не на самих элементарных теориях, а на их фрагментах, соответствующих префиксным классам $\Sigma_k$ и $\Pi_k$, где $k$ --- натуральное число; с этой целью А.\ Нис  слегка модифицирует метод интерпретаций. Поскольку можно считать, что $\Sigma_k \cup \Pi_k \subseteq \Sigma_{k+1} \cap \Pi_{k+1}$, естественно возникает задача поиска минимальных неразрешимых (максимальных разрешимых) фрагментов теорий. В частности, А.\ Нис доказывает наследственную неразрешимость $\Sigma_2$-фраг\-мен\-та теории конечных простых графов. Ввиду разрешимости $\Pi_2$-фраг\-мен\-та, это говорит о том, что $\Sigma_2$ --- минимальный неразрешимый фрагментом рассматриваемой теории. Стоит отметить, что трансляции, которые использовались в более ранних работах по элементарным теориям, зачастую далеки от оптимальных с точки зрения префиксной сложности. Поэтому в \cite{Nies} пришлось придумывать новые, более оптимальные трансляции. Больше информации можно найти в \cite{Nies} и \cite{Speranski-2016-JLC}.

Цель настоящей работы --- достижение прогресса в поиске минимальных неразрешимых фрагментов базовых элементарных теорий. Будет доказано, что класс всех конечных двудольных графов является $\Sigma_1$-ин\-тер\-пре\-ти\-ру\-е\-мым в классе всех пар эквивалентностей на общем конечном носителе; отсюда будет получаться наследственная неразрешимость $\Sigma_2$-те\-о\-рии второго класса. Далее, будет построена $\Sigma_1$-ин\-тер\-пре\-та\-ция класса всех пар эквивалентностей на общем конечном носителе в классе всех пар, состоящих из эквивалентности и линейного порядка на общем конечном носителе; это даст наследственную неразрешимость $\Sigma_2$-те\-о\-рии второго класса. Соответствующие результаты являются в известном смысле оптимальными ввиду разрешимости $\Pi_2$-те\-о\-рий рассматриваемых классов. 


\section{Предварительные сведения}

Напоминаем, что \emph{сигнатурой} называется набор предикатных, функциональных и константных символов с указанием местности для предикатных и функциональных.

Для сигнатуры $\sigma$ обозначим через $K^{\sigma}$ класс всех $\sigma$-струк\-тур. Здесь и далее мы ограничимся рассмотрением сигнатур без функциональных символов.

\begin{defi}
Для сигнатур $\sigma_1$ и $\sigma_2$ под \emph{$\sigma_2$-схемой в $\sigma_1$} понимается набор $\sigma_2$-формул:
\begin{itemize}

\item $\Phi_U \left( x, \overline{y} \right)$;

\item $\Phi_R \left( x_1, \dots, x_n, \overline{y} \right)$ и $\Phi_{\neg R} \left( x_1, \dots, x_n, \overline{y} \right)$ для каждого предикатного символа $R$ из $\sigma_1$, где $n$ равно местности $R$.

\end{itemize}
Здесь $\overline{y}$ --- некоторый кортеж переменных.
\end{defi}

\begin{defi} \label{1.1}
Пусть $\mathcal{K}_1 \subseteq K_{\sigma_1}$ и $\mathcal{K}_2 \subseteq K_{\sigma_2}$. Говорят, что $\mathcal{K}_1$ \emph{интерпретируем с па\-ра\-мет\-рами} (сокращённо \emph{и.п.}) в $\mathcal{K}_2$, если существует $\sigma_2$-схема в $\sigma_1$ такая, что для любой $\mathfrak{A} \in \mathcal{K}_1$ найдутся $\mathfrak{B} \in \mathcal{K}_2$ 
и кортеж $\overline{p}$ из $B$, удовлетворяющие следующим условиям:
\begin{enumerate}

\item $B' := \left\{ b \mid b \in B\ \text{и}\ \mathfrak{B} \vDash {\Phi_U \left( b, \overline{p} \right)} \right\}$ не~пусто;

\item для любых $n$-мест\-но\-го предикатного символа $R$ из $\sigma_1$ и $\left( b_1, \dots, b_n \right) \in {\left( B' \right)}^n$,
\[
{\mathfrak{B} \vDash \Phi_{\neg R} \left( b_1, \dots, b_n, \overline{p} \right)}
\quad \Longleftrightarrow \quad
{\mathfrak{B} \nvDash \Phi_R \left( b_1, \dots, b_n, \overline{p} \right)} ;
\]

\item $\mathfrak{A}$ изоморфно $\sigma_1$-структуре $\mathfrak{B}'$ с носителем $B'$, где каждый $n$-мест\-ный предикатный символ~$R$ из $\sigma_1$ интерпретируется по правилу
\[
{\mathfrak{B}' \vDash {R \left( b_1, \dots, b_n, \overline{p} \right)}} \quad :\Longleftrightarrow \quad {\mathfrak{B} \vDash {\Phi_R \left( b_1,
\dots, b_n, \overline{p} \right)}} .
\]

\end{enumerate}
Разумеется, в случае, когда кортеж $\overline{y}$ в схеме пуст, слова <<с параметрами>> опускают.\footnote{В
русскоязычной литературе вместо <<интерпретируем>> часто используется <<относительно элементарно определим>>; см.\ \cite[глава 5]{ershov1980}.}
\end{defi}

\begin{rem}
Приведённое выше определение не является наиболее общим, хотя для наших целей его будет достаточно. В частности, 
в схемы нередко добавляют специальную $\sigma_2$-формулу $\Phi_{=} \left( x_1, x_2, \overline{y} \right)$. Тогда дополнительно требуется, чтобы множество пар
\[
F\ :=\ \left\{ \left( b_1, b_2 \right) \in B^2 \mid \mathfrak{B} \vDash {\Phi_= \left( b_1, b_2, \overline{p} \right)} \right\}
\]
было отношением конгруэнции на $\mathfrak{B}'$, а условие <<$\mathfrak{A}$ изоморфно $\mathfrak{B}'$>> заменяется на несколько более мягкое <<$\mathfrak{A}$ изоморфно ${\mathfrak{B}'} /_{F}$>>. Кроме того, нередко используют многомерные интерпретации, когда элементы $A$ моделируются с помощью кортежей элементов $B$.
\end{rem}

Для $\mathcal{K} \subseteq K^{\sigma}$ обозначим через $\mathrm{Th} \left( \mathcal{K} \right)$ совокупность всех $\sigma$-предложений, истинных во всех структурах из $\mathcal{K}$. Кроме того, для каждого $k \in \mathbb{N}$ положим
\begin{align*}
{\Sigma_k \text{-} \mathrm{Th} \left( \mathcal{K} \right)}\
&:=\
{\mathrm{Th} \left( \mathcal{K} \right) \cap \Sigma_k} ,\\
{\Pi_k \text{-} \mathrm{Th} \left(\mathcal{K} \right)}\
&:=\
{\mathrm{Th} \left( \mathcal{K} \right) \cap \Pi_k} .
\end{align*}
Тут $\Sigma_k$ обозначает формулы вида ${\exists \vec{x}_1}\, {\forall \vec{x}_2}\, {\exists \vec{x}_3}\, \ldots\, \vec{x}_k\, \Phi$, где префикс содержит ровно $k$ блоков однотипных кванторов, а $\Phi$ --- бескванторная формула. Аналогично для $\Pi_k$, но меняя местами $\exists$ и $\forall$; таким образом, $\Pi_k$-фор\-му\-лы можно воспринимать как отрицания $\Sigma_k$-фор\-мул.

Говорят, что класс $\mathcal{K}_1$ является \emph{$\Sigma_k$-интерпретируемым (с параметрами)} в классе $\mathcal{K}_2$, если все формулы в схеме из определения \ref{1.1} лежат в $\Sigma_k$; см.\ \cite{Nies}. Аналогично для $\Pi_k$.

В дальнейшем мы будем всегда подразумевать некоторую эффективную (гёделевскую) нумерацию для формул данной сигнатуры и, как следствие, отождествлять формулы с их номерами.

\begin{defi}
Для сигнатуры $\sigma$ обозначим через $\mathrm{Val}_{\sigma}$ множество всех общезначимых $\sigma$-пред\-ло\-же\-ний. Говорят, что множество $\sigma$-формул $\Gamma$ \emph{наследственно неразрешимо} (сокращённо \emph{н.н.}), если для любого множества $\sigma$-формул $\Delta$,
\[
{\mathrm{Val}_{\sigma} \cap \Gamma}\ \subseteq\ \Delta\ \subseteq\ \Gamma
\quad \Longrightarrow \quad
\Delta \enskip \text{неразрешимо}
\]
--- это фактически означает, что $\mathrm{Val}_{\sigma} \cap \Gamma$ и дополнение $\Gamma$ вычислимо неотделимы.
\end{defi}

\begin{lemm}[{о трансляции; cм.\ \cite[лемма~3.1]{Nies}}] \label{l1.1}
Для любых чисел $r$, $l$ и классов $\mathcal{K}_1$, $\mathcal{K}_2$:
\begin{align*} \label{l1.1}
{\left. {\begin{array}{r}
\mathcal{K}_1 ~ \text{является} ~ \Sigma_k \text{-и.\,п.\ в} ~ \mathcal{K}_2,\\
{\Pi_{r + 1} \text{-} {\mathrm{Th} \left( \mathcal{K}_1 \right)}} ~ \text{является н.\,н.}\\
\end{array}} \right\}}
\quad &\Longrightarrow \quad {\Pi_{r + k} \text{-} {\mathrm{Th} \left( \mathcal{K}_2 \right)}} ~ \text{является н.\,н.} ;\\[0.5em]
{\left. {\begin{array}{r}
\mathcal{K}_1 ~ \text{является} ~ \Sigma_k \text{-и.\ в} ~ \mathcal{K}_2,\\
{\Sigma_r \text{-} {\mathrm{Th} \left( \mathcal{K}_1 \right)}} ~ \text{является н.\,н.} \\
\end{array}} \right\}}
\quad &\Longrightarrow \quad {\Sigma_{r + k - 1} \text{-} {\mathrm{Th} \left( \mathcal{K}_2 \right)}} ~ \text{является н.\,н.}
\end{align*}
\end{lemm}

Для удобства обозначим через $K^{\sigma}_\mathrm{fin}$ класс всех конечных $\sigma$-струк\-тур. Отметим один достаточно~простой факт:

\begin{prop}[{см.\ \cite[раздел 6.2.2]{Boerger&Graedel&Gurevich-1997}, например}]
Предположим, что $\sigma$ не~содержит функциональных символов. Тогда $\Pi_2 \text{-} \mathrm{Th} \left( K^{\sigma} \right)$ разрешима и совпадает с $\Pi_2 \text{-} \mathrm{Th} \left( K^{\sigma}_\mathrm{fin} \right)$.
\end{prop}

Отсюда можно легко получить следующее.

\begin{prop} \label{p2}
Предположим, что $\sigma$ не~содержит функциональных символов. Пусть $\mathcal{K}$ --- произвольный класс $\sigma$-cтрук\-тур, аксиоматизируемый конечным числом $\Sigma_2$-пред\-ло\-же\-ний. Тогда $\Pi_2 \text{-} \mathrm{Th} \left( \mathcal{K} \right)$ разрешима и совпадает с $\Pi_2 \text{-} \mathrm{Th} \left( \mathcal{K} \cap K^{\sigma}_\mathrm{fin} \right)$.
\end{prop}

\begin{myproof}
Ясно, что $\mathcal{K}$ можно аксиоматизировать посредством одного $\Sigma_2$-пред\-ло\-же\-ния $\Theta$. Тогда для любого $\Pi_2$-пред\-ло\-же\-ния $\Phi$,
\[
\Phi\ \in\ {\mathrm{Th} \left( \mathcal{K} \right)}
\quad \Longleftrightarrow \quad
{\Theta \rightarrow \Phi}\ \in\ {\mathrm{Th} \left( K^{\sigma} \right)} . \tag{$\dag$}
\]
При этом $\Theta \rightarrow \Phi$, очевидно, логически эквивалентно $\Pi_2$-пред\-ло\-же\-нию. Значит, $\Pi_2 \text{-} \mathrm{Th} \left( \mathcal{K} \right)$ сводится к $\Pi_2 \text{-} \mathrm{Th} \left( K^{\sigma} \right)$, а потому разрешима. Совпадение $\Pi_2 \text{-} \mathrm{Th} \left( \mathcal{K} \right)$ и $\Pi_2 \text{-} \mathrm{Th} \left( \mathcal{K} \cap K^{\sigma}_\mathrm{fin} \right)$ вытекает из совпадения $\Pi_2 \text{-} \mathrm{Th} \left( K^{\sigma} \right)$ и $\Pi_2 \text{-} \mathrm{Th} \left( K^{\sigma}_\mathrm{fin} \right)$ по модулю ($\dag$).
\end{myproof}

\begin{defi}
Рассмотрим сигнатуру $\sigma := \{E^2\}$. Под {\it простым графом} мы понимаем $\sigma$-струк\-ту\-ру, удовлетворяющую предложению
\[
{\forall x}\, {\neg E (x, x)} \wedge
{\forall x}\, {\forall y}\, {\left( E (x, y) \leftrightarrow E (y, x) \right)} .
\]
Обозначим через ${G}_\mathrm{fin}$ класс всех конечных простых графов.
\end{defi}

Наследственная неразрешимость теории конечных простых графов была доказана в \cite{lavrov63}.\footnote{Под
<<эффективной неотделимостью>> в названии статьи в действительности понимается <<вычислимая неотделимость>> в современной терминологии; ср.\ \cite{Speranski-2016-JLC}.}
Позже в \cite[лемма~4.2]{Nies} была доказана наследственная неразрешимость $\Sigma_2 \text{-} \mathrm{Th} \left( {G}_\mathrm{fin} \right)$. 

С помощью леммы \ref{l1.1} в различных работах (например, \cite{Nies} и \cite{Speranski-2016-JLC}) из наследственной неразрешимости теории конечных простых графов выводится наследственная неразрешимость многих других теорий.

\begin{defi}
Рассмотрим сигнатуру $\sigma := \{L^1, R^1, E^2\}$. Под {\it двудольным графом} понимают $\sigma$-структуру, удовлетворяющую предложению    
\[
{\forall x}\, {\left( L (x) \leftrightarrow {\neg R (x)} \right)} \wedge
{\forall x}\, {\forall y}\, {\left( E (x, y) \rightarrow {L (x) \wedge R (y)} \right)} .
\]
В этом случае носитель каждой модели можно естественным образом разбить на \emph{левую} и \emph{правую доли}. Обозначим через ${BiG}_\mathrm{fin}$ класс всех конечных двудольных графов.
\end{defi}

Наследственная неразрешимость $\Sigma_2 \text{-} \mathrm{Th} \left( {BiG}_\mathrm{fin} \right)$ была доказана в \cite[лемма~4.5]{Nies}. 

\begin{rem}
Из доказательства Ниса следует, что наследственно неразрешимым также будет $\Sigma_2$-фраг\-мент теории конечных двудольных графов, содержащих не менее трёх вершин в каждой из долей. Для упрощения рассуждений, мы будем использовать именно этот класс, обозначив его за ${BiG}^{\ast}_\mathrm{fin}$.
\end{rem}

\begin{defi}
Рассмотрим сигнатуру $\sigma := \{P^2, Q^2\}$. Под {\it моделью двух эквивалентностей} понимают $\sigma$-структуру, удовлетворяющую предложению
\begin{multline*}
{\forall x}\, {P (x, x)} \wedge {\forall x}\, {\forall y}\, {\left( P (x, y) \rightarrow P (y, x) \right)} \wedge {\forall x}\, {\forall y}\, {\forall z}\, {\left( {P (x, y) \wedge P (y, z)} \rightarrow P(x, z) \right)} \mathop{\wedge}\\
{\forall x}\, {Q (x, x)} \wedge {\forall x}\, {\forall y}\, {\left( Q (x, y) \rightarrow Q (y, x) \right)} \wedge {\forall x}\, {\forall y}\, {\forall z}\, {\left( {Q (x, y) \wedge Q (y, z)} \rightarrow Q(x, z) \right)} .
\end{multline*}
Обозначим через ${2Eq}_\mathrm{fin}$ класс всех конечных моделей двух эквивалентностей, через ${2Eq}$ --- класс всех моделей двух эквивалентностей.
\end{defi}

В доказательстве теоремы 3.6 в \cite{Speranski-2016-JLC} было отмечено, что из \cite[cc.\ 273--274]{ershov1980} можно получить наследственную неразрешимость $\Sigma_3 \text{-} \mathrm{Th} \left( {2Eq}_\mathrm{fin} \right)$. Вместе с тем утверждение \ref{p2} гарантирует разрешимость $\Pi_2 \text{-} \mathrm{Th} \left( {2Eq}_\mathrm{fin} \right).$

\begin{defi}
Рассмотрим сигнатуру $\sigma := \left\{ <^2, \approx^2 \right\}$. Под \emph{моделью линейного порядка и эквивалентности} понимают $\sigma$-структуру, в которой $<$ интерпретируется как линейный порядок, а $\approx$ --- отношение эквивалентности.
Обозначим через ${LEq}_\mathrm{fin}$ класс всех конечных моделей линейного порядка и эквивалентности, через ${LEq}$ --- класс всех моделей линейного порядка и эквивалентности.
\end{defi}

Как показывает доказательство предложения 8.6.10 в \cite{ershov2011}, ${2Eq}_\mathrm{fin}$ будет $\Sigma_2$-ин\-тер\-пре\-ти\-ру\-е\-мым в ${LEq}_\mathrm{fin}$. Отсюда с помощью леммы \ref{l1.1} получается наследственная неразрешимость $\Sigma_4 \text{-} \mathrm{Th} \left( {LEq}_\mathrm{fin} \right)$. Вместе с тем утверждение \ref{p2} гарантирует разрешимость $\Pi_2 \text{-} \mathrm{Th} \left( {LEq}_\mathrm{fin} \right)$.


\section{О теории двух эквивалентностей}

В настоящем разделе мы сначала построим $\Sigma_1$-ин\-тер\-пре\-та\-цию с параметрами класса всех конечных двудольных графов в классе всех конечных моделей двух эквивалентностей. Затем эта конструкция будет специальным образом модифицирована, чтобы избавиться от параметров. В итоге мы докажем, что $\Sigma_2 \text{-} \mathrm{Th} \left( {2Eq}_\mathrm{fin} \right)$ является минимальным неразрешимым фрагментом рассматриваемой теории.

\begin{lemm} \label{2EqP3}
${BiG}^{\ast}_\mathrm{fin}$ является $\Sigma_1$-ин\-тер\-пре\-ти\-ру\-е\-мым с параметрами в ${2Eq}_\mathrm{fin}$.
\end{lemm}

\begin{myproof}
Пусть $\mathfrak{A}$ --- конечный двудольный граф с не менее чем тремя вершинами в каждой из долей. Для удобства занумеруем его левую и правую доли:
\[
L^{\mathfrak{A}}\ =\ {\{ l_1, \ldots, l_m \}}
\quad \text{и} \quad
R^{\mathfrak{A}}\ =\ {\{ r_1, \ldots, r_n \}} .
\]
Заметим, что $A$ (т.е.\ носитель $\mathfrak{A}$) совпадает с $L^{\mathfrak{A}} \sqcup R^{\mathfrak{A}}$.

Далее, зафиксируем специальное множество
\[
C\ := \{ c_P, c_N, c_L, c_R \} ;
\]
его элементы будут играть роль параметров в дальнейшем. Кроме того, заведём три дополнительных множества, каждое из которых содержит ровно $m \cdot n$ элементов:
\begin{align*}
S_L\ &:=\ \left\{ s_L^{ij} \mid 1 \le i \le m,\ 1 \le j \le n \right\} ;\\
S_R\ &:=\ \left\{ s_R^{ij} \mid 1 \le i \le m,\ 1 \le j \le n \right\} ;\\
S_E\ &:=\ \left\{ s_E^{ij} \mid 1 \le i \le m,\ 1 \le j \le n \right\} .
\end{align*}
Разумеется, мы предполагаем, что $L^{\mathfrak{A}}$, $R^{\mathfrak{A}}$, $C$, $S_L$, $S_R$ и $S_E$ попарно не~пересекаются.

Построим по $\mathfrak{A}$ конечную модель двух эквивалентностей $\mathfrak{B}$ с носителем
\[
B\ :=\ {L^{\mathfrak{A}} \cup R^{\mathfrak{A}} \cup S_L \cup S_R \cup S_E \cup C} ,
\]
в которой отношения эквивалентности $P^{\mathfrak{B}}$ и $Q^{\mathfrak{B}}$ определяются посредством разбиений $B$ на классы следующим образом.

Классы эквивалентности по $P^{\mathfrak{B}}$:
\begin{itemize}

\item $\{ c_L, l_1, l_2, \ldots, l_m \}$;

\item $\{ c_R, r_1, r_2, \ldots, r_n \}$;

\item $\{ s_L^{ij}, s_R^{ij}, s_E^{ij} \}$ для любых $i \in \{ 1, \ldots, m \}$ и $j \in \{ 1, \ldots, n \}$;

\item $\{ c_P \}$;

\item $\{c_N \}$.

\end{itemize}

Классы эквивалентности по $Q^{\mathfrak{B}}$:
\begin{itemize}

\item $\{ l_i \} \cup \{ s_L^{ij} \mid 1 \le j \le n \}$ для каждого $i \in \{ 1, \ldots, m \}$;

\item $\{ r_j \} \cup \{ s_R^{ij} \mid 1 \le i \le m \}$ для каждого $j \in \{ 1, \ldots, n \}$;

\item $\{ c_P \} \cup \{ s_E^{ij} \mid \mathfrak{A} \vDash E \left( l_i, r_j \right) \}$ для любых $i \in \{ 1, \ldots, m \}$ и $j \in \{ 1, \ldots, n \}$;

\item $\{ c_N \} \cup \{ s_E^{ij} \mid \mathfrak{A} \vDash \neg E \left( l_i, r_j \right) \}$ для любых $1 \le i \le m$, $1 \le j \le n$;

\item $\{ c_L \}$;

\item $\{ c_R \}$.

\end{itemize}

Например, для двудольного графа с долями $\{ l_1, l_2 \}$ и $\{ r_1, r_2 \}$ и рёбрами $\left( l_1, r_1 \right)$ и $\left( l_2, r_2 \right)$ разбиения на классы можно схематично представить следующим образом:\footnote{Тут
для наглядности мы ограничиваемся примером с двумя вершинами в каждой из долей.}

\begin{figure}[!htbp]
  \centering
  \begin{minipage}[b]{0.4\textwidth}
    \includegraphics[width=\textwidth]{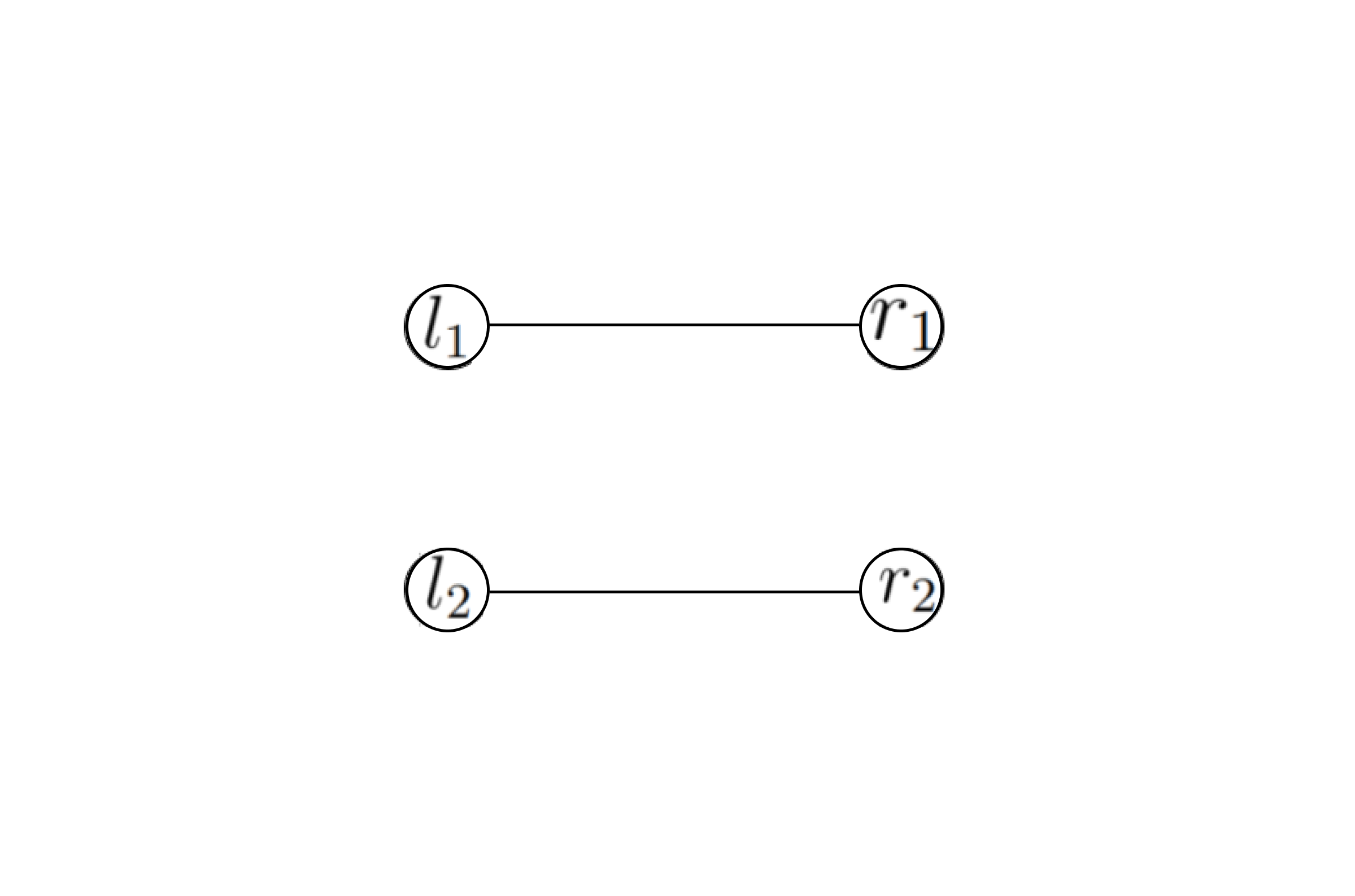}
    \caption{Исходный граф.}
  \end{minipage}
  \hfill
  \begin{minipage}[b]{0.4\textwidth}
    \includegraphics[width=\textwidth]{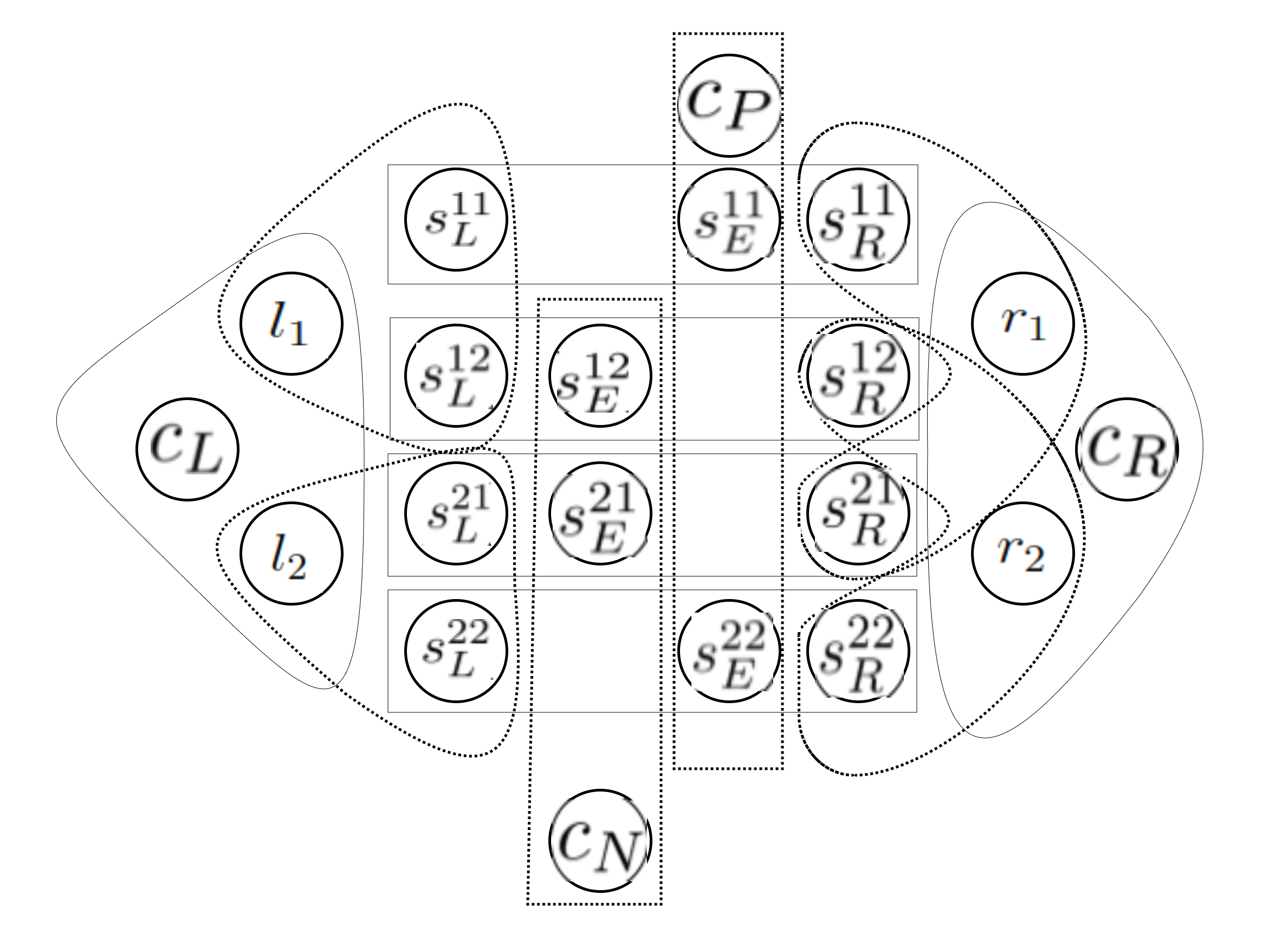}
    \caption{Два разбиения $B$ на классы.}
  \end{minipage}
\end{figure}

С целью упрощения записи введём сокращения
\begin{align*}
{\overline{P} \left( x_1, \ldots, x_k \right)}\ :=\ {\bigwedge_{1 \le i, j \le k} P \left( x_i, x_j \right)} ,\\
{\overline{Q} \left( x_1, \ldots, x_k \right)}\ :=\ {\bigwedge_{1 \le i, j \le k} Q \left( x_i, x_j \right)} .
\end{align*}
Теперь опишем схему трансляции из сигнатуры $\left\{ L^1, R^1, E^2 \right\}$ в сигнатуру $\left\{ P^2, Q^2 \right\}$ с параметрами $\overline{y} = \left( y_L, y_R, y_P, y_N \right)$:
\begin{align*}
{\Phi_U \left( x, \overline{y} \right)}\ := \enskip
&{\left( \neg Q \left( x, y_L \right) \wedge P \left( x, y_L \right) \right) \vee \left( \neg Q \left( x, y_R \right) \wedge P \left( x, y_R \right) \right)} ;\\[0.5em]
{\Phi_L \left( x, \overline{y} \right)}\ := \enskip
&{P \left( x, y_L \right)} ;\\[0.5em]
{\Phi_{\neg L} \left( x, \overline{y} \right)}\ := \enskip
&{\neg P \left( x, y_L \right)} ;\\[0.5em]
{\Phi_{R} \left( x, \overline{y} \right)}\ := \enskip
&{P \left( x, y_R \right)} ;\\[0.5em]
{\Phi_{\neg R} \left( x, \overline{y} \right)}\ := \enskip
&{\neg P \left( x, y_R \right)} ;\\[0.5em]
{\Phi_{E} \left( x_L, x_R, \overline{y} \right)}\ := \enskip
&\Phi_L \left( x_L \right) \wedge \Phi_R \left( x_R \right) \mathop{\wedge}\\
&{\exists u_L, u_R, u_E}\, {\left( Q \left( x_L, u_L \right) \wedge Q \left( x_R, u_R \right) \wedge \overline{P} \left( u_L, u_R, u_E \right) \wedge Q \left( u_E, y_P \right) \right)} ;\\[0.5em]
{\Phi_{\neg E} \left( x_L, x_R, \overline{y} \right)}\ := \enskip
&\Phi_{\neg L} \left( x_L \right) \vee \Phi_{\neg R} \left( x_R \right) \mathop{\vee}\\
&{\exists u_L, u_R, u_E}\, {\left( Q \left( x_L, u_L \right) \wedge Q \left( x_R, u_R \right) \wedge \overline{P} \left( u_L, u_R, u_E \right) \wedge Q \left( u_E, y_N \right) \right)} .
\end{align*}
Нетрудно убедиться, что все условия из определения \ref{1.1} окажутся выполнены, если для исходной $\mathfrak{A}$ взять $\mathfrak{B}$, построенную выше, с элементами $c_L$, $c_R$, $c_P$, $c_N$ в качестве значений для параметров $y_L$, $y_R$, $y_P$, $y_N$.
Таким образом, мы показали $\Sigma_1$-ин\-тер\-пре\-ти\-ру\-е\-мость с параметрами ${BiG}_\mathrm{fin}$ в ${2Eq}_\mathrm{fin}$.
\end{myproof}

\begin{theo}
${BiG}_\mathrm{fin}^{\ast}$ является $\Sigma_1$-ин\-тер\-пре\-ти\-ру\-е\-мым (без параметров) в ${2Eq}_\mathrm{fin}$.
\end{theo}

\begin{myproof}
Ниже будет приведена трансляция, которая в целом напоминает построенную ранее, но уже без параметров. Основная идея состоит в том, чтобы вместо $c_L$, $c_R$, $c_P$ и $c_N$ (из предыдущего доказательства) использовать специальные $\Sigma_1$-оп\-ре\-де\-ли\-мые множества; при этом $C$ придётся заменить на более сложное множество.

Пусть $\mathfrak{A}$ --- конечный двудольный граф с не менее чем тремя вершинами в каждой из долей. Пусть $L^{\mathfrak{A}}$ и $R^{\mathfrak{A}}$ занумерованы так же, как и раньше. Определения $S_L$, $S_R$ и $S_E$ также не~меняются. Вместо $C$ введём четыре специальных множества:
\begin{align*}
C_L\ &:=\
\{ {c}_L^{1}, \ldots, {c}_L^{4} \} \cup \{ \widetilde{c}_L^{\ 1}, \ldots, \widetilde{c}_L^{\ 4} \} ;\\
C_R\ &:=\
\{ {c}_R^{1}, \ldots, {c}_R^{\ 3} \} \cup
\{ \widetilde{c}_R^{\ 1} \ldots \widetilde{c}_R^{\ 5} \} ;\\
C_P\ &:=\
\{ c_P^{1}, \ldots, c_P^{6} \} \cup \{ \widetilde{c}_P^{\ 1}, \ldots, \widetilde{c}_P^{\ 6} \} ;\\
C_N\ &:=\
\{ c_N^{1}, \ldots, c_N^{7} \} \cup \{ \widetilde{c}_N^{\ 1}, \ldots, \widetilde{c}_N^{\ 7} \} .
\end{align*}

Построим по $\mathfrak{A}$ конечную модель двух эквивалентностей $\mathfrak{B}$ с носителем
\[
B\ :=\
L^{\mathfrak{A}} \cup R^{\mathfrak{A}} \cup S_L \cup S_R \cup S_E \cup C_L \cup C_R \cup C_P \cup C_N ,
\]
в которой отношения эквивалентности $P^{\mathfrak{B}}$ и $Q^{\mathfrak{B}}$ определяются посредством разбиений $B$ на классы следующим образом.

Классы эквивалентности по $P^{\mathfrak{B}}$:
\begin{itemize}

\item $\{ c_L^1, \ldots, c_L^4, l_1, l_2, \ldots, l_m \}$;

\item $\{ \widetilde{c}_L^{\ 1}, \ldots, \widetilde{c}_L^{\ 4} \}$;

\item $\{ c_R^1, \ldots, c_R^3, r_1, r_2, \ldots, r_n\}$;

\item $\{ \widetilde{c}_R^{\ 1}, \ldots, \widetilde{c}_R^{\ 5} \}$;

\item $\{s_L^{ij}, s_R^{ij}, s_E^{ij}\}$ для любых $i \in \{ 1, \ldots, m \}$ и $j \in \{ 1, \ldots, n \}$;

\item $\{ c_P^1, \ldots, c_P^6 \}$;

\item $\{ \widetilde{c}_P^{\ 1}, \ldots, \widetilde{c}_P^{\ 6} \}$;

\item $\{ c_N^{1}, \ldots, c_N^7 \}$;

\item $\{ \widetilde{c}_N^{\ 1}, \ldots, \widetilde{c}_N^{\ 7} \}$.\footnote{Здесь $m$ и $n$ суть количества элементов в левой и правой долях $\mathfrak{A}$ соответственно.}

\end{itemize}

Классы эквивалентности по $Q^{\mathfrak{B}}$:
\begin{itemize}

\item $\{ l_i \} \cup \{s_L^{ij} \mid 1 \le j \le n\}$ для каждого $i \in \{ 1, \ldots, m \}$;

\item  $\{ r_j \} \cup \{s_R^{ij} \mid 1 \le i \le m\}$ для каждого $j \in \{ 1, \ldots, n \}$;

\item $\{ c_P^{6} \} \cup \{ s_E^{ij} \mid \mathfrak{A} \vDash E(l_i, r_j)\}$ для любых $i \in \{ 1, \ldots, m \}$ и $j \in \{ 1, \ldots, n \}$;

\item $\{ c_N^{7} \} \cup \{ s_E^{ij} \mid \mathfrak{A} \vDash \neg E(l_i, r_j)\}$ для любых $i \in \{ 1, \ldots, m \}$ и $j \in \{ 1, \ldots, n \}$;

\item $\{ c_L^k, \widetilde{c}_L^{\ k} \}$ для всех $k \in \{ 1, 2, 3, 4 \}$;

\item $\{ c_R^k, \widetilde{c}_R^{\ k} \}$ для всех $k \in \{ 1, 2, 3 \}$, а также синглетоны из неиспользованных элементов $C_R$;

\item $\{ c_P^k, \widetilde{c}_P^{\ k} \}$ для всех $k \in \{ 1, 2 \}$, а также синглетоны из неиспользованных элементов $C_P$;

\item $\{c_N^1, \widetilde{c}_N^{\ 1} \}$, а также синглетоны из неиспользованных элементов $C_N$.

\end{itemize}
Заметим, что для любых $b_1, b_2 \in B$, если $b_1 \ne b_2$, то $\left( b_1, b_2 \right) \not \in P^{\mathfrak{B}}$ или $\left( b_1, b_2 \right) \not \in Q^{\mathfrak{B}}$. Стало быть, равенство в $\mathfrak{B}$ можно определить посредством формулы $P \left( x, y \right) \wedge Q \left( x, y \right)$.

С целью упрощения записи введём дополнительные сокращения: для любых чисел $k_1$, $k_2$ и $k_3$ таких, что  $k_3 \le \min \left\{ k_1, k_2 \right\}$, возьмём
\begin{align*}
{\Theta_{k_1, k_2, k_3} \left( \vec{x}, \vec{y} \right)}\ := \enskip
&\overline{P} \left( x_1, \ldots, x_{k_1} \right) \wedge \overline{P} \left( y_1, \ldots, y_{k_2} \right) \wedge {\neg P \left( x_1, y_1 \right)} \mathop{\wedge}\\
&\bigwedge_{1 \le i, j \le k_1} {\neg Q \left( x_i, x_j \right)} \wedge \bigwedge_{1 \le i, j \le k_2} {\neg Q \left( y_i, y_j \right)} \mathop{\wedge}\\
&\bigwedge_{1 \le i \le k_3}  {Q \left( x_i, y_i \right)} ,
\end{align*}
где $\vec{x}$ и $\vec{y}$ обозначают $\left( x_1, \ldots, x_{k_1} \right)$ и $\left( y_1, \ldots, y_{k_2} \right)$ соответственно. На самом деле, нам будет нужно лишь конечное число таких формул.

%
%
%

Пусть $B_1$ и $B_2$ --- различные классы эквивалентности по $P^{\mathfrak{B}}$. Мы будем говорить, что \emph{$B_1$ и $B_2$ порождают $k$ связей по $Q^{\mathfrak{B}}$}, если $\left( B_1 \times B_2 \right) \cap Q^{\mathfrak{B}}$ содержит в точности $k$ элементов. В частности:
\begin{itemize}

\item $\{c_L^1 \ldots c_L^4, l_1, l_2, \ldots l_m\}$ и $\{\widetilde{c}_L^{\ 1}, \ldots \widetilde{c}_L^{\ 4}\}$ порождают четыре связи по $Q^{\mathfrak{B}}$;

\item $\{c_R^1 \ldots c_R^3, r_1, r_2, \ldots r_n\}$ и $\{\widetilde{c}_R^{\ 1}, \ldots \widetilde{c}_R^{\ 5}\}$ порождают три связи по $Q^{\mathfrak{B}}$;

\item $\{c_P^1, \ldots, c_P^6\}$ и $\{\widetilde{c}_P^{\ 1}, \ldots, \widetilde{c}_P^{\ 6}\}$ порождают две связи по $Q^{\mathfrak{B}}$;

\item $\{c_N^1 \ldots c_N^7\}$ и $\{\widetilde{c}_N^{\ 1}, \ldots \widetilde{c}_N^{\ 7}\}$ порождают одну связь по $Q^{\mathfrak{B}}$.

\end{itemize}
Нетрудно проверить, что для всякого означивания $\nu$ в $\mathfrak{B}$:
\begin{itemize} 

\item если $\mathfrak{B} \vDash \Theta_{5, 4, 4} \left( \vec{x}, \vec{y} \right) \left[ \nu \right]$ (т.е.\ $\Theta_{5, 4, 4} \left( \vec{x}, \vec{y} \right)$ истинна в $\mathfrak{B}$ при означивании $\nu$), то
\[
{\{ c_L^1, \ldots, c_L^4 \}}\ =\
{\{ {\nu \left( x_1 \right)}, \ldots, {\nu \left( x_4 \right)} \}}
\quad \text{и} \quad
{\{ \widetilde{c}_L^1, \ldots, \widetilde{c}_L^{4} \}}\ =\
{\{ {\nu \left( y_1 \right)}, \ldots, {\nu \left( y_4 \right)} \}} ;
\]

\item если $\mathfrak{B} \vDash \Theta_{6, 5, 3} \left( \vec{x}, \vec{y} \right) \left[ \nu \right]$, то
\[
{\{ c_R^1, \ldots, c_R^3 \}}\ =\ {\{ {\nu \left( x_1 \right)}, \ldots, {\nu \left( x_3 \right)} \}}
\quad \text{и} \quad
{\{ \widetilde{c}_R^1, \ldots, \widetilde{c}_R^{5} \}}\ =\ {\{ {\nu \left( y_1 \right)}, \ldots, {\nu \left( y_5 \right)} \}} ;
\]

\item если $\mathfrak{B} \vDash \Theta_{6, 6, 2} \left( \vec{x}, \vec{y} \right) \left[ \nu \right]$, то
\[
{\{ c_P^1, \ldots, c_P^6 \} \sqcup \{ \widetilde{c}_P^{\ 1}, \ldots, \widetilde{c}_P^{\ 6} \}}\ =\
{\{ {\nu \left( x_1 \right)}, \ldots, {\nu \left( x_6 \right)}, {\nu \left( y_1 \right)}, \ldots, {\nu \left( y_6 \right)} \}} ;
\]

\item если $\mathfrak{B} \vDash \Theta_{7, 7, 1} \left( \vec{x}, \vec{y} \right) \left[ \nu \right]$, то
\[
{\{ c_N^1, \ldots, c_N^7 \} \sqcup \{ \widetilde{c}_N^{\ 1}, \ldots, \widetilde{c}_N^{\ 7} \}}\ =\
{\{ {\nu \left( x_1 \right)}, \ldots, {\nu \left( x_7 \right)}, {\nu \left( y_1 \right)}, \ldots, {\nu \left( y_7 \right)} \}} .
\]
       
\end{itemize}
Далее, рассмотрим формулы:
\begin{align*}
{\Psi_L \left( x \right)}\ := \enskip
&{\exists y_1, \ldots, y_4}\, {\exists z_1, \ldots, z_4}\, {\left( \Theta_{5, 4, 4} \left( y_1, \ldots, y_4, x; z_1, \ldots, z_4 \right) \right)} ;\\[0.5em]
{\Psi_R \left( x \right)}\ := \enskip
&{\exists y_1, \ldots, y_5}\, {\exists z_1, \ldots, z_5}\, {\left( \Theta_{6, 5, 3} \left( y_1, \ldots, y_5, x; z_1, \ldots, z_5  \right) \right)} ;\\[0.5em]
{\Psi_P \left( x \right)}\ := \enskip
&{\exists y_1, \ldots, y_5}\, {\exists z_1, \ldots, z_6}\, {\left( \Theta_{6, 6, 2} \left( x, y_1, \ldots, y_5; z_1, \ldots, z_6 \right) \mathop{\lor} \right.}\\
&\hspace{15.1em} {\left. \Theta_{6, 6, 2} \left( y_1, \ldots, y_5, x; z_1, \ldots, z_6 \right) \right)} ;\\[0.5em]
{\Psi_N \left( x \right)}\ := \enskip
&{\exists y_1, \ldots, y_6}\, {\exists z_1, \ldots, z_7}\, {\left( \Theta_{7, 7, 1} \left( x, y_1, \ldots, y_6; z_1, \ldots, z_7 \right) \mathop{\lor} \right.}\\
&\hspace{15.1em} {\left. \Theta_{7, 7, 1} \left( y_1, \ldots, y_6, x; z_1, \ldots, z_7 \right) \right)} .
\end{align*}
Нетрудно убедиться, что множества $\{ l_1, \ldots, l_m \}$, $\{ r_1, \ldots, r_n \}$, $\{ c_P^1, \ldots, c_P^6 \} \cup \{ \widetilde{c}_P^{\ 1}, \ldots, \widetilde{c}_P^{\ 6} \}$ и\linebreak $\{ C_N^1, \ldots, c_N^7 \} \cup \{ \widetilde{c}_N^{\ 1}, \ldots, \widetilde{c}_N^{\ 7} \}$ будут $\Sigma_1$-оп\-ре\-де\-ли\-мы в $\mathfrak{B}$ посредством формул $\Psi_L \left( x \right)$, $\Psi_R \left( x \right)$, $\Psi_P \left( x \right)$ и $\Psi_N \left( x \right)$ соответственно. Наконец, мы готовы описать искомую схему трансляции без параметров:
\begin{align*}
{\Phi_L \left( x \right)}\ := \enskip
&{\Psi_L \left( x \right)} ;\\[0.5em]
{\Phi_{R} \left( x \right)}\ := \enskip
&{\Psi_R \left( x \right)} ;\\[0.5em]
{\Phi_{\neg L} \left( x \right)}\ := \enskip
&{\Psi_R \left( x \right)} ;\\[0.5em]
{\Phi_{\neg R} \left( x \right)}\ := \enskip
&{\Psi_L \left( x \right)} ;\\[0.5em]
{\Phi_U \left( x \right)}\ := \enskip
&{\Phi_L \left( x \right) \vee \Phi_R \left( x \right)} ;\\[0.5em]
{\Phi_E \left( x, y \right)}\ := \enskip
&{\Phi_L \left( x \right)} \wedge {\Phi_R \left( y \right)} \mathop{\wedge}\\
&{\exists u_L, u_R, u_E, u_P}\, {\left( \Psi_P \left( u_P \right) \wedge Q \left( x, u_L \right) \wedge Q \left( y, u_R \right) \wedge \overline{P} \left( u_L, u_R, u_E \right) \wedge Q \left( u_E, u_P \right) \right)} ;\\[0.5em]
{\Phi_{\neg E} \left( x, y \right)}\ := \enskip
&\Phi_{\neg L}(x) \vee \Phi_{\neg R}(y) \mathop{\vee}\\
&{\exists u_L, u_R, u_E, u_N}\, {\left( \Psi_N \left( u_N \right) \wedge Q \left( x, u_L \right) \wedge Q \left( y, u_R \right) \wedge \overline{P} \left( u_L, u_R, u_E \right) \wedge Q \left(u_E, u_N \right) \right)} .
\end{align*}
Ясно, что все эти формулы эквивалентны $\Sigma_1$-фор\-му\-лам. Более того, можно проверить, что соответствующая $\mathfrak{B}'$ (см.\ определение \ref{1.1}) будет изоморфна $\mathfrak{A}$. В итоге мы показали $\Sigma_1$-ин\-тер\-пре\-ти\-ру\-е\-мость ${BiG}_\mathrm{fin}$ в ${2Eq}_\mathrm{fin}$ без параметров.
\end{myproof}

\begin{coro}
$\Sigma_2 \text{-} \mathrm{Th} \left( {2Eq}_\mathrm{fin} \right)$ наследственно неразрешима.
\end{coro}

\begin{myproof}
Получается из наследственной неразрешимости $\Sigma_2 \text{-} \mathrm{Th} \left( {BiG}^{\ast}_\mathrm{fin} \right)$ и $\Sigma_1$-ин\-тер\-пре\-ти\-ру\-е\-мо\-сти ${BiG}^{\ast}_\mathrm{fin}$ в ${2Eq}_\mathrm{fin}$ с помощью леммы \ref{l1.1}.
\end{myproof}

\begin{rem}
В приведённой трансляции не использовались элементы, попарно равные друг-другу по обеим эквивалентностям. Далее при доказательстве сводимости теории двух эквивалентностей к другим теориям мы будем считать, что любые два различных элемента можно отличить по одной из эквивалентностей.
\end{rem}

%



%



%



%

Наконец, с помощью более оптимальной трансляции можно улучшить некоторые из оценок в \cite[Теорема~3.6]{Speranski-2016-JLC}:

\begin{coro}
$\Pi_5$-фрагмент теории свободных дистрибутивных решёток с конечным числом порождающих наследственно неразрешим.
\end{coro}

\begin{myproof}
В \cite[cc.\ 279--281]{ershov1980} показано, что ${2Eq}_\mathrm{fin}$ является $\Sigma_3$-интерпретируемым с параметрами в рассматриваемом классе решёток. Далее применяем лемму \ref{l1.1}.
\end{myproof}

\begin{coro}
$\Pi_5$-фрагмент теории конечных групп перестановок наследственно неразрешим.
\end{coro}

\begin{myproof}
В \cite[cc.\ 283--285]{ershov1980} показано, что ${2Eq}_\mathrm{fin}$ является $\Sigma_3$-интерпретируемым с параметрами в рассматриваемом классе групп. Далее применяем лемму \ref{l1.1}.
\end{myproof}


\section{О теории линейного порядка и эквивалентности}

Как было показано в \cite[cc.\ 339--340]{ershov2011}, класс всех конечных моделей двух эквивалентностей является $\Sigma_2$-ин\-тер\-пре\-ти\-ру\-е\-мым с параметрами в классе всех конечных моделей линейного порядка и эквивалентности; отсюда можно получить наследственную неразрешимость $\Pi_4 \text{-} \mathrm{Th} \left( {LEq}_{\mathrm{fin}} \right)$, но эта оценка по-прежнему не является оптимальной. В настоящем разделе мы сначала построим $\Sigma_1$-ин\-тер\-пре\-та\-цию с параметрами ${2Eq}_{\mathrm{fin}}$ в ${LEq}_{\mathrm{fin}}$. Затем эта конструкция будет модифицирована, чтобы избавиться от параметров. В итоге мы докажем, что $\Sigma_2 \text{-} \mathrm{Th} \left( {LEq}_\mathrm{fin} \right)$ является минимальным неразрешимым фрагментом рассматриваемой теории. Стоит отметить, что часть идей для построения трансляции с параметрами взята из доказательства предложения~8.6.10 в \cite{ershov2011}.

\begin{lemm} \label{leqp}
${2Eq}_\mathrm{fin}$ является $\Sigma_1$-ин\-тер\-пре\-ти\-ру\-е\-мым с параметрами в ${LEq}_\mathrm{fin}$.
\end{lemm}

\begin{myproof}
Пусть $\mathfrak{A}$ --- конечная модель двух эквивалентностей. Занумеруем классы эквивалентности, на которые разбивают множество $A$ отношения $P^{\mathfrak{A}}$ и $Q^{\mathfrak{A}}$: 
\begin{align*}
A\ =\ {A_1^P \sqcup \ldots \sqcup A_n^P}
\quad \text{и} \quad
A\ =\ {A_1^Q \sqcup \ldots \sqcup A_m^Q} ,
\end{align*}
где $n$ и $m$ --- количества классов эквивалентности по $P^\mathfrak{A}$ и $Q^\mathfrak{A}$ соответственно. Возьмём
\[
S\ :=\
{\{ s_k^i \mid 1 \le i \le n,\ 1 \le k \le 3 \}} .
\]
Кроме того, зафиксируем некоторый новый элемент $a^{\ast}$.

Построим по $\mathfrak{A}$ конечную модель линейного порядка и эквивалентности $\mathfrak{B}$ с носителем
\[
B\ :=\ {A \cup S \cup \{ a^{\ast}\}} .
\]
следующим образом.

Классы эквивалентности по $\approx^\mathfrak{B}$:
\begin{itemize}

\item $\{ s_0^i, s_1^i, s_2^i \}$ для каждого $i \in \left\{ 1, \ldots, n \right\}$;

\item $A_j^Q$ для каждого $j \in \left\{ 1, \ldots, m \right\}$;

\item $\left\{ a^{\ast} \right\}$.

\end{itemize}
Данное разбиение соответствует определению $\approx^\mathfrak{B}$.

Теперь определим $<^{\mathfrak{B}}$. Мы будем считать, что каждый класс $A_i^P$ линейно упорядочен произвольным образом. Тогда $<^\mathfrak{B}$ зададим так:
\[
s_0^1\ <\ s_0^2\ <\ \ldots\ <\ s_0^n\ <\ a^{\ast}\ <\ s_1^1\ <\ A^P_1\ <\ s_2^1\ <\ s_1^2\ <\ A^P_2\ <\ s_2^2\ <\ \ldots\ < s_1^n\ <\ A^P_n\ <\ s_2^n ,
\]
где запись $s^i_k < A^P_i$ означает, что $s^i_k$ меньше всех элементов $ A^P_i$.

Ниже приведён пример того, как выглядят линейный порядок и эквивалентность в $\mathfrak{B}$ для конкретной $\mathfrak{A}$.

\begin{figure}[!htbp]
  \centering
  \begin{minipage}[b]{0.45\textwidth}
    \centering
    \includegraphics[width=50mm]{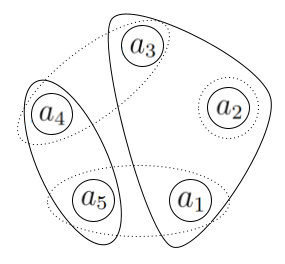}
    \caption{Два разбиения $A$ на классы. Здесь сплошная линия отвечает за эквивалентность по $P^{\mathfrak{A}}$, а пунктирная --- по $Q^{\mathfrak{A}}$.}
  \end{minipage}
  \hfill
  \begin{minipage}[b]{0.50\textwidth}
    \includegraphics[width=\textwidth]{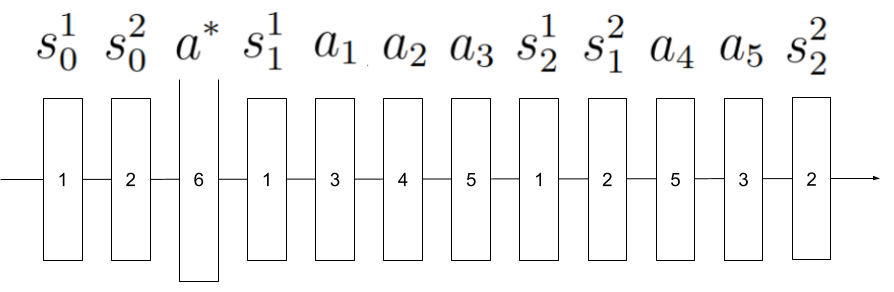}
    \caption{Линейный порядок и эквивалентность на $B$. Элементы расположены от меньшего к большему в соответствии с $<^{\mathfrak{B}}$; два элемента эквивалентны по $\approx^{\mathfrak{B}}$, если и только если им отвечает одно и то же число.}
  \end{minipage}
\end{figure}

Теперь опишем схему трансляции из сигнатуры $\left\{ P^2, Q^2 \right\}$ в сигнатуру $\left\{ \approx^2, <^2 \right\}$ с параметром $y$:
\begin{align*}
{\Phi_U \left( x, y \right)}\ := \enskip
&{\exists z_0, z_1, z_2}\, {\left( z_0 \approx z_1 \approx z_2 \wedge z_0 < y < z_1 < x < z_2 \right)} ;\\[0.5em]
{\Phi_{\neg U} \left( x, y \right)}\ := \enskip
&{\exists z_0}\, {\left( z_0 \approx x \wedge \neg (z_0 > y) \right)} ;\\[0.5em]
{\Phi_Q \left( x_1, x_2, y \right)}\ := \enskip
&{x_1 \approx x_2} ;\\[0.5em]
{\Phi_{\neg Q} \left( x_1, x_2, y \right)}\ := \enskip
&{\neg {x_1 \approx x_2}} ;\\[0.5em]
{\Phi_P \left( x_1, x_2, y \right)}\ := \enskip
&{\exists z_0, z_1, z_2}\, {\left( z_0 \approx z_1 \approx z_2 \wedge z_0 < y < z_1 < x_1 < z_2 \wedge z_1 < x_2 < z_2 \right)} ;\\[0.5em]
{\Phi_{\neg P} \left( x_1, x_2, y \right)}\ := \enskip
&{\exists z}\, {\left( \Phi_{\neg U} \left( z, y \right) \wedge \left( x_1 < z < x_2 \vee x_2 < z < x_1 \right) \right)} .
\end{align*}
Нетрудно убедиться, что все условия из определения \ref{1.1} окажутся выполнены, если для исходной $\mathfrak{A}$ взять $\mathfrak{B}$, построенную выше, с элементом $a^{\ast}$ в качестве значения для параметра $y$. В итоге мы показали $\Sigma_1$-ин\-тер\-пре\-ти\-ру\-е\-мость с параметрами ${2Eq}_\mathrm{fin}$ в ${LEq}_\mathrm{fin}$.
\end{myproof}

\begin{rem}
В доказательстве последней леммы элемент $a^{\ast}$ можно определить посредством конъюнкции двух $\Pi_1$-фор\-му\-л:
\[
{\forall x}\, {\left( x > y \rightarrow {\neg x \approx y} \right)} \wedge {\forall x_1}\, {\forall x_2}\, {\left( {x_1 < x_2 < y} \rightarrow {\neg {x_1 \approx x_2}} \right) } .
\]
Отсюда следует $\Sigma_2$-ин\-тер\-пре\-ти\-ру\-е\-мость ${2Eq}_\mathrm{fin}$ в ${LEq}_\mathrm{fin}$.
\end{rem}

\begin{theo}
${2Eq}_\mathrm{fin}$ является $\Sigma_1$-ин\-тер\-пре\-ти\-ру\-е\-мым (без параметров) в ${LEq}_\mathrm{fin}$.
\end{theo}

\begin{myproof}
Наша цель --- избавиться от параметра в построенной ранее трансляции. Основная идея состоит в том, чтобы новая модель линейного порядка и эквивалентности давала ровно один класс эквивалентности, содержащий более трёх элементов, который и будет играть роль $\Sigma_1$-оп\-ре\-де\-ли\-мо\-го параметра.

Пусть $\mathfrak{A}$ --- конечная модель двух эквивалентностей. Мы будем считать, что
\[
A\ =\ {\left\{ a_1, \ldots, a_{|A|} \right\}} .
\]
Пусть классы эквивалентности по $P^\mathfrak{A}$ и $Q^\mathfrak{A}$ занумерованы так же, как и раньше; $n$ и $m$ --- количества классов эквивалентности по $P^\mathfrak{A}$ и $Q^\mathfrak{A}$ соответственно. Вместо $a^{\ast}$ введём специальное множество
\[
C\ :=\ {\left\{ c_k \mid 1 \le k \le 4 \right\}} .
\]
В дополнение к $S$ введём
\[
R\ :=\ {\left\{ r_k^j \mid 1 \le j \le m,\ 1 \le k \le 3 \right\}} .
\]
Кроме того, нам понадобится занумерованная копия $A$:
\[
\widetilde{A}\ =\ {\left\{ \widetilde{a}_1, \ldots, \widetilde{a}_{|A|} \right\}} .
\]
Разумеется, на $\widetilde{A}$ можно перенести два исходных разбиения:
\[
\widetilde{A}\ =\ {\widetilde{A}_1^{P} \sqcup \ldots \sqcup \widetilde{A}_n^{P}}
\quad \text{и} \quad
\widetilde{A}\ =\ {\widetilde{A}_1^{Q} \sqcup \ldots \sqcup \widetilde{A}_m^{Q}} .
\]
    
Построим по $\mathfrak{A}$ конечную модель линейного порядка и эквивалентности $\mathfrak{B}$ с носителем
\[\
B\ :=\
{A \cup \widetilde{A} \cup S \cup R \cup C}
\]
следующим образом.

Классы эквивалентности по $\approx^\mathfrak{B}$:
\begin{itemize}

\item $\{ a_l, \widetilde{a}_l \}$ для всех $l \in \left\{ 1, \ldots, {\left| A \right|} \right\}$; 

\item $\{ s_0^i, s_1^i, s_2^i \}$ для каждого $i \in \left\{ 1, \ldots, n \right\}$;

\item$\{ r_0^j, r_1^j, r_2^j \}$ для каждого $j \in \left\{ 1, \ldots, m \right\}$;

\item $\left\{ c_1, c_2, c_3, c_4 \right\}$.

\end{itemize}
Данное разбиение соответствует определению $\approx^\mathfrak{B}$.

Линейный порядок $<^{\mathfrak{B}}$ зададим так: \pagebreak
\begin{gather*}
s^1_0\ <\ \ldots\ <\ s^n_0\ < r^1_0\ <\ \ldots\ <\ r^m_0\ <\ c_1\ <\ c_2\ <\ c_3\ <\ c_4\ <\\[0.5em]
<\ [ r^1_1 < A^Q_1 < r^1_2 ]\ <\ \ldots\ <\  [ r^m_1 < A^Q_m < r^m_2 ]\ <\ [ s^1_1 < \widetilde{A}^{P}_1 < s^1_2 ]\ <\ \ldots\ <\ [ s_1^n < \widetilde{A}^{P}_n < s_2^n ].
\end{gather*}
Здесь мы считаем, что каждый класс $\widetilde{A}^{P}_i$ и каждый класс $A_j^Q$ линейно упорядочены произвольным образом. Квадратные скобки добавлены для удобства чтения.

Ниже приведён пример того, как выглядят линейный порядок и эквивалентность в $\mathfrak{B}$ для конкретной $\mathfrak{A}$; ср.\ с предыдущим доказательством.

\begin{figure}[!htbp]
  \centering
  \begin{minipage}[b]{0.45\textwidth}
    \centering
    \includegraphics[width=50mm]{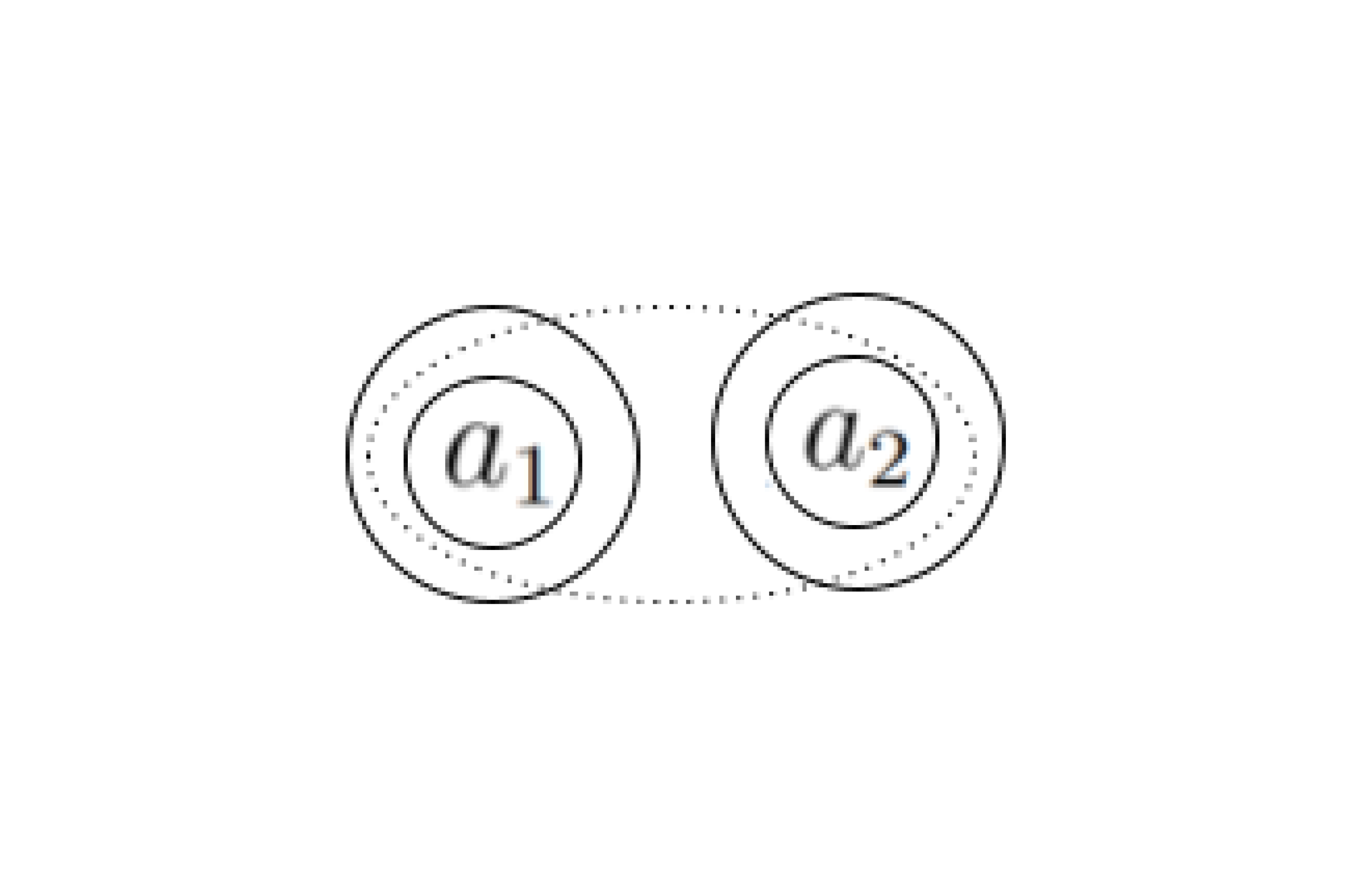}
    \caption{Два разбиения $A$ на классы. Здесь сплошная линия отвечает за эквивалентность по $Q^{\mathfrak{A}}$, а пунктирная --- по $P^{\mathfrak{A}}$.}
  \end{minipage}
  \hfill
  \begin{minipage}[b]{0.50\textwidth}
    \includegraphics[width=\textwidth]{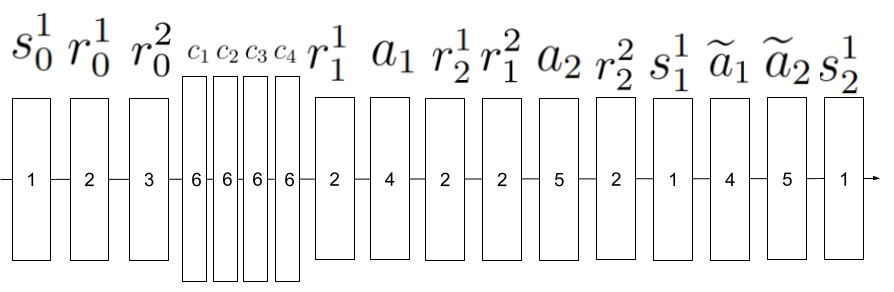}
    \caption{Линейный порядок и эквивалентность на $B$.}
  \end{minipage}
\end{figure}

Так как ровно один класс эквивалентности по $\approx^{\mathfrak{B}}$ содержит более трёх элементов, мы можем определить $c_1$ в $\mathfrak{B}$ посредством $\Sigma_1$-фор\-му\-лы
\[
{\Psi \left( x \right)}\ :=\
{\exists x_1}\, {\exists x_2}\, {\exists x_3}\, {\left( x \approx x_1 \approx x_2 \approx x_3 \wedge x < x_1 < x_2 < x_3 \right)} .
\]
Мы уже готовы описать искомую схему трансляции без параметров:
\begin{align*}
{\Phi_U \left( x \right)}\ :=\
&{\exists y, y_0, y_1, y_2, z}\, {\left( \Psi \left( y \right) \wedge {y_0 < y < y_1 < x < y_2} \mathop{\wedge} \right.}\\
&{\left. {y_0 \approx y_1 \approx y_2} \wedge {x \approx z} \wedge {x < z} \right)} ;\\[0.5em]
{\Phi_P \left( x_1, x_2 \right)} :=\
&{\exists y, y_0, y_1, y_2, z_1, z_2}\, ( \Psi \left( y \right) \wedge {y_0 < y < y_1 < y_2} \mathop{\wedge}\\
&{y_0 \approx y_1 \approx y_2} \wedge {z_1 \approx x_1} \wedge {z_2 \approx x_2} \wedge {y_1 < z_1 < y_2} \wedge {y_1 < z_2 < y_2} ) ; \\[0.5em]
{\Phi_{\neg P} \left( x_1, x_2 \right)} :=\
&{\exists y, y_0, y_1, y_2, z_1, z_2}\, ( \Psi \left( y \right) \wedge {y_0 < y < y_1 < y_2} \mathop{\wedge}\\
&{y_0 \approx y_1 \approx y_2} \wedge {z_1 \approx x_1} \wedge {z_2 \approx x_2} \wedge \left( {z_1 < y_2 < z_2} \vee {z_2 < y_2 < z_1} \right) ) ;\\[0.5em]
{\Phi_Q \left( x_1, x_2 \right)} :=\
&{\exists y, y_0, y_1, y_2}\, ( \Psi \left( y \right) \wedge {y_0 < y < y_1 < y_2} \mathop{\wedge}\\
&{y_0 \approx y_1 \approx y_2} \wedge {y_1 < x_1 < y_2} \wedge {y_1 < x_2 < y_2} ) ;\\[0.5em]
{\Phi_{\neg Q} \left( x_1, x_2 \right)} :=\
&{\exists y, y_0, y_1, y_2}\, ( \Psi \left( y \right) \wedge {y_0 < y < y_1 < y_2} \mathop{\wedge}\\
&{y_0 \approx y_1 \approx y_2} \wedge \left( {x_1 < y_2 < x_2} \vee {x_2 < y_2 < x_1} \right) ) .
\end{align*}
Ясно, что все эти формулы эквивалентны $\Sigma_1$-фор\-му\-лам. Более того, можно проверить, что соответствующая $\mathfrak{B}'$ (см.\ определение \ref{1.1}) будет изоморфна $\mathfrak{A}$. Отсюда вытекает $\Sigma_1$-ин\-тер\-пре\-ти\-ру\-е\-мость ${2Eq}_\mathrm{fin}$ в ${LEq}_\mathrm{fin}$ без параметров.
\end{myproof}

\begin{coro} \label{2-LEQ}
$\Sigma_2 \text{-} \mathrm{Th} \left( {LEq}_\mathrm{fin} \right)$ наследственно неразрешима.
\end{coro}

\begin{myproof}
Получается из наследственной неразрешимости $\Sigma_2 \text{-} \mathrm{Th} \left( {2Eq}_\mathrm{fin} \right)$ и $\Sigma_1$-ин\-тер\-пре\-ти\-ру\-е\-мо\-сти ${2Eq}_\mathrm{fin}$ в ${LEq}_\mathrm{fin}$ с помощью леммы \ref{l1.1}.
\end{myproof}

%



%


%






\section{Заключение}

В настоящей работе были найдены минимальные неразрешимые префиксные фрагменты теорий классов ${2Eq}_{\mathrm{fin}}$ и ${LEq}_{\mathrm{fin}}$. Поскольку ${2Eq}_{\mathrm{fin}}$ часто используется для доказательства наследственной неразрешимости теорий других классов, это позволяет улучшить ряд более ранних результатов о неразрешимости префиксных фрагментов. Тем не менее, улучшенные оценки часто не~являются оптимальными и нуждаются в уточнении.

В частности, естественным образом возникает идея рассмотрения класса  ${2L}_\mathrm{fin}$, состоящего из всех пар линейных порядков на общем конечном носителе. Ясно, что $\Pi_2$-фраг\-мент соответствующей теории будет разрешим. С другой стороны, её $\Pi_4$-фрагмент является наследственно неразрешимым --- это можно получить из доказательства предложения 8.6.11 в \cite{ershov2011}, где приведена $\Sigma_2$-ин\-тер\-пре\-та\-ция ${LEq}_{\mathrm{fin}}$ в ${2L}_{\mathrm{fin}}$ с параметрами, и нашего следствия \ref{2-LEQ}. На самом деле, от параметров можно избавиться, что даст неразрешимость $\Sigma_3$-фраг\-мен\-та. Тем не менее, открытым остаётся вопрос о разрешимости $\Sigma_2$- и $\Pi_3$-фраг\-мен\-тов.




\bigskip
В заключение автор выражает глубокую благодарность своему научному руководителю Станиславу Олеговичу Сперанскому за проявленное внимание к работе, чрезвычайно ценные замечания и неизменно плодотворные обсуждения материала.


\addcontentsline{toc}{section}{Список литературы}

\newpage
\begin{flushleft}

\end{flushleft}

\end{document}